\documentstyle[11pt]{article}
\begin{document}
\setlength{\baselineskip}{17pt}
\begin{center}
{\Large {\bf \mbox{} \\ \mbox{} \\ New Travelling Wave Solutions for the \\
Fisher-KPP Equation with General Exponents}}
\end{center}

\mbox{}

\begin{center}
{\large {\sc Ariel S\'{a}nchez--Vald\'{e}s \hspace{1cm}
Benito Hern\'{a}ndez--Bermejo \footnote{\normalsize Corresponding author. E-mail: {\tt bhernandez@escet.urjc.es } \\ \mbox{} \hspace{0.35cm} Telephone: (+34) 91 488 73 91. Fax: (+34) 91 488 73 38.}  }} \\
\mbox{} \\
Universidad Rey Juan Carlos. \\
E.S.C.E.T. Campus de M\'{o}stoles. Edificio Departamental II. \\
Calle Tulip\'{a}n S/N. 28933--M\'{o}stoles--Madrid. Spain. \\
\end{center}

\mbox{}

\mbox{}

\mbox{}

\noindent{\bf Abstract---}
The Fisher-KPP equation with general nonlinear diffusion and arbitrary kinetic orders in the reaction terms is considered. The existence of oscillatory travelling wave solutions is proved
for this model. Conditions for the existence of such solutions are provided.

\mbox{}

\mbox{}

\mbox{}

\noindent {\bf Keywords---} Fisher-KPP equation, reaction-diffusion, travelling wave solutions.

\vfill

\pagebreak
The following reaction-diffusion equation is to be considered in what follows:
\begin{equation}
u_t= \kappa (u^{m-1}u_x)_{x}+ \alpha u^p- \beta u^q \; , \qquad x \in I \! \! R \; , \; \; t>0
\label{ec0}
\end{equation}
with $\alpha , \beta , \kappa , m>0$ and $p,q \in I \! \! R$. After the change of variables $x \to ax$, $t\to bt$, $u\to lu$, where $a=( \kappa l^{m-p}/\alpha)^{1/2}$, $b=(l^{1-p}/\alpha)$, and $l=(\beta/\alpha)^{1/(p-q)}$, equation (\ref{ec0}) becomes:
\begin{equation}
u_t=(u^{m-1}u_x)_{x}+ u^p- u^q \; , \qquad x \in I \! \! R \; , \; \; t>0
\label{ec1}
\end{equation}
The existence of travelling wave solutions of (\ref{ec1}) will be investigated in this work. Recall that these are solutions of the form $u(x,t)=f(x-ct)=f(\xi)$, where $c\in I \! \! R$ is the wave velocity, $\xi\in I \! \! R $. Such solutions satisfy
\[
f(\xi)\geq 0 \: , \quad \lim_{\xi\to\infty} f(\xi)=0 \: , \quad \lim_{\xi\to -\infty} f(\xi)=1
\]
and consequently will be termed right travelling wave solutions (RTW from now on).

Many particular cases of (\ref{ec1}) have been considered in detail in the literature. For
instance, equation (\ref{ec1}) with $m=1$, $p=1$, and $q=2$ was proposed by Fisher [\ref{Fi}]
as a model for the evolution of a gene favourable for a given population. Later, Kolmogorov {\em et al.} [\ref{KPP}] demonstrated that there exists a threshold value $c^*=2$ for the Fisher equation such that RTW exist if and only if the wave velocity $c$ satisfies $c \geq c^*$.

Additionally, for the more general situation $p\leq q$, $m>0$ and $m+p>0$ it is well-known that there is also a constant $c^*$ (of value not yet explicitly determined in general) such that RTW exist if and only if $c\geq c^*$ and $m+p\geq 2$. Specifically, this result was demonstrated for the homogeneous diffusion case $m=1$ by Aronson and Weinberger [\ref{AW}], while the remaining non-homogenous possibilities were analysed by several authors [\ref{GK},\ref{B},\ref{PV}]. For the additional case of equations with general diffusion and reaction coefficients, see also [\ref{GK},\ref{MM1},\ref{SM}]. It is worth recalling in this context that the consideration of general  reaction exponents $p$ and $q$ such as those used in
(\ref{ec1}) is of central interest in the framework of applied biochemical modelling (see
[\ref{VO}--\ref{M}] and references therein for a detailed revision).

In the complementary case $p>q$, a relevant consequence of the work due to Gilding and Kersner [\ref{GK}] is that if $\: (m+q)>0$, then there exists a value $c^*<0$ such that RTW are present for $f \in [0,1]$ provided $c \leq c^*$. The main result of the present work is that the previous solutions do not exhaust all the possibilities for that situation:

\mbox{}

\noindent{\sc Theorem 1.} {\sl If $p>q$ and $(m+q)>0$ in equation (\ref{ec1}), there exists a RTW if and only if $c<0$. Moreover, if $\mid c \mid \geq 2\sqrt{p-q}=\mid c^* \mid $ then $f \in [0,1]$, while if $0< \mid c \mid < \mid c^* \mid $ then the RTW oscillates around the limit value 1 with strictly decreasing amplitude when $\xi \rightarrow - \infty$.}

\mbox{}

\noindent{\sc Proof.} After transformation $u(x,t)=f(\xi)$ equation (\ref{ec1}) becomes
\begin{equation}
(f^{m-1}f')'+cf'+f^p-f^q=0
\label{ec2}
\end{equation}
where $f(\xi)\to 1$ if $\xi\to-\infty$ and $f(\xi)\to 0$ if $\xi\to\infty$. Equation
(\ref{ec2}) is to be interpreted in the weak sense, namely $f$ and $(f^{m-1}f')$ are continuous in $ I \! \! R $, and the equation is satisfied in its standard integral version.

\mbox{}

\noindent{\sc Remark 1.} {\sl Note that if $c$ and $\xi$ are replaced by $-c$ and $-\xi$ respectively, equation (\ref{ec2}) remains unaltered. Therefore looking for RTW (that evolve from $1$ to $0$) with negative velocity is equivalent to investigating travelling wave solutions with positive velocity from $0$ to $1$.}

\mbox{}

\noindent Let us now convert equation (\ref{ec2}) to a first order system whose phase plane is to be investigated. For this, two cases must be distinguished:

\mbox{}

\noindent {\em I. Case $\: (m+q)>2$:} Defining
\[
X=f^{m+q-2} \: , \:\:\:\: Y=f^{m-2}f' \: , \:\:\:\: X^{\frac{m-1}{m+q-2}}d\tau=d\xi
\]
the following system arises:
\begin{equation}
\label{eq3} \left\{
\begin{array}{ll}
X'=\gamma XY=P(X,Y)\\ [4pt] Y'=-Y(Y+c)+X-X^k=Q(X,Y)
\end{array}
\right.
\end{equation}
where $k=\frac{m+p-2}{m+q-2}>1$ and $\gamma =m+q-2$. Taking Remark 1 into account, let us
suppose $c>0$ in order to look for trajectories leaving $X=0$ and entering $X=1$ (RTW with positive velocity), or trajectories leaving $X=1$ and going to $X=0$ (RTW with negative  velocity). Note also that only trajectories belonging to the invariant region $S=\{(X,Y):X>0\}$ are relevant. In addition, a simple application of Dulac's criterion with weight function
$B=X^{\frac{2}{\gamma}-1}$ in $S$ shows that there are no periodic orbits in $S$.
This fact and the Poincar\'{e}-Bendixon theorem are to be considered in order to investigate system (\ref{eq3}).

\noindent Now the local analysis of (\ref{eq3}) leads to the following fixed points: {\em (a)\/} $P_0=(0,0)$ is a saddle-node point with saddle region belonging to $S$; {\em (b)\/} $P_1=(0,-c)$ is a saddle; and finally {\em (c)\/} $P_2=(1,0)$ which is a stable node if $c\geq 2\sqrt{p-q}=c^*$, and a focus if $c<c^*$. These results together with the analysis of the dynamics at the fixed points $(0, \infty)$ and $(0,- \infty)$ show that the only admissible trajectories are those joining $P_0$ and $P_2$, which implies that there are no RTW with positive velocity.

\noindent Next we demonstrate that if $c>0$ there exists a unique trajectory joining $P_2$ and $P_0$. Let us call ${\Gamma_0}(X,Y)$ the trajectory leaving $P_0$ and ${\Gamma_1}(X,Y)$ the trajectory leaving $P_1$. Then a standard analysis of (\ref{eq3}) shows that the following statements $C_1$, $C_2$, $C_3$ and $C_4$ are true in $S$:
\begin{description}
  \item[$C_1:$] $X'>0$ if $Y>0$ and $X'<0$ if $Y<0$.
  \item[$C_2:$] ${\Gamma_0}(X,Y)$ intersects the $X$ axis. Let $X_0$ be the first of such
    intersections. Then $X_0\geq 1$.
  \item[$C_3:$] ${\Gamma_1}(X,Y)$ does not intersect the interval $[0,X_0)$.
  \item[$C_4:$] Either ${\Gamma_1}(X,Y)$ intersects the straight line $Y=0$ or the $X$
    component of ${\Gamma_1}(X,Y)$ tends to $\infty$ when $\tau\to -\infty$.
\end{description}
\noindent In addition, $C_1$, $C_2$, $C_3$ and $C_4$ together imply that either ${\Gamma_0}(X,Y)\to P_2$ or ${\Gamma_0}(X,Y)={\Gamma_1}(X,Y)$. We next have:

\mbox{}

\noindent{\sc Lemma 1.} {\sl ${\Gamma_0}(X,Y)$ enters $P_2$ if $c>0$.}

\mbox{}

\noindent{\sc Proof.} Consider the equation for the trajectories $Y(c,X)$ of system
(\ref{eq3}):
\begin{equation}
\frac{dY}{dX}=\frac{-Y(Y+c)+X-X^k}{\gamma XY}
\label{ec-tray}
\end{equation}
From now on $Y_{\Gamma}(c,X)$ will denote the solution of (\ref{ec-tray}) corresponding to the solution curve $\Gamma$. Then the fact that $Y_{\Gamma_0}(c,X)$ behaves like $(X/c)$ in the neighbourhood of $(0,0)$ implies that $c_2>c_1$ and then
$Y_{\Gamma_0}(c_1,X)>Y_{\Gamma_0}(c_2,X)$ in the vicinity of $X=0$. Let us now prove that
$Y_{\Gamma_0}(c_1,X)$ and $Y_{\Gamma_0}(c_2,X)$ do not intersect for $Y>0$. For this suppose the opposite, and let $(X_1,Y_1)$ be the intersection point. This implies that
$\frac{d}{dX}Y_{\Gamma_0}(c_1,X_1) \leq \frac{d}{dX}Y_{\Gamma_0}(c_2,X_1)$ and then $c_1 \geq c_2$, which is a contradiction. This shows that if $X_{01}$ and $X_{02}$ are the points
where $Y_{\Gamma_0}(c_1,X)$ and $Y_{\Gamma_0}(c_2,X)$ intersect the $X$-axis for the first
time, respectively, then $X_{01}\geq X_{02}$. Similarly it can be proved that if $X_{11}$ and $X_{12}$ are the points where $Y_{\Gamma_1}(c_1,X)$ and $Y_{\Gamma_1}(c_2,X)$ intersect the
$X$-axis for the first time, with $c_2>c_1$, then $X_{11}\leq X_{12}$. In addition we have that if $c=0$ then the trajectories are explicit and equal to $Y_{\Gamma_1}(0,X)$, and they satisfy
the equation $Y^2-\frac{2X}{2+\gamma}+\frac{2X^k}{2+\gamma k}=0$. Then the intersection point with the $X$-axis here is $X_0=(\frac{2+\gamma k}{2+\gamma})^{\frac{1}{k-1}}>1$. Therefore if
$c_2>c_1$ then $X_{12}\geq X_{11}\geq X_0\geq X_{01}\geq X_{02}$. To conclude we only need to
prove that if $c_2>c_1$ and $X_{01}>1$, then $X_{01}\ne X_{02}$. Suppose the contrary,
$X_{01}= X_{02}$. We know that $Y_{\Gamma_0}(c_1,X)>Y_{\Gamma_0}(c_2,X)$ if $X<X_{01}$. To
simplify, let us rename $Y_{\Gamma_0}(c_1,X)$ and $Y_{\Gamma_0}(c_2,X)$, $Y_1$ and $Y_2$ respectively. We then have
\[
\frac{d}{dX}(Y_1-Y_2) =\frac{Y_1Y_2((c_2-c_1)-(Y_1-Y_2))+
X(X^{k-1}-1)(Y_1-Y_2)}{\gamma X Y_1Y_2}>0
\]
if $X$ is close enough to $X_{01}>1$, because $(Y_1-Y_2)\to 0$ if $X\to X_{01}$. Consequently
we arrive to a contradiction regarding that $Y_{\Gamma_0}(c_1,X)$ and $Y_{\Gamma_0}(c_2,X)$ do coincide at $X_{01}$. This completes the proof of Lemma 1. \hfill {\sc Q.E.D.}

\mbox{}

\noindent Let us prove the second part of Case I, namely that $X_0=1$ if $c\geq2\sqrt{p-q}$. We only need to show that $X_0=1$ if $c=2\sqrt{p-q}$. For this, we prove that the trajectory does not leave the region:
\[
G=\{(X,Y):X\geq0,Y\geq0,Y+a(X-1)\leq0\} \: , \:\:\: a=\frac{c}{2(m+q-2)}
\]
Due to condition $C_2$, ${\Gamma_0}(X,Y)$ does not intersect the coordinate axes in $G$, and then we only need to prove that it does not intersect the line $\{ r:Y+a(X-1)=0 \}$. For this we shall prove that $\mbox{\bf n}\cdot\mbox{\bf V}<0$, where $\mbox{\bf n}=(a,1)$ is the normal vector of $r$, and $\mbox{\bf V}$ is the restriction of the vector field (\ref{eq3}) to $r$. Thus:
\[
\mbox{\bf n}\cdot\mbox{\bf V}=a(m+q-2)XY-Y(Y+c)-X^k+X
\]
Since $Y=a(1-X)$ we have:
\[
\mbox{\bf n}\cdot\mbox{\bf V}=-X^2a^2(m+q-1)+X(a^2(m+q)+ca+1)-a^2-ca-X^k \equiv R(X)
\]
Note that $R(0)<0$ and $R(1)=0$, and if $R'(X)\geq 0$ in $[0,1]$, then $R(X)\leq 0$ which demonstrates the result, and the proof that $R'(X)\geq 0$ is:
\[
R'(X)=-2a^2(m+q-1)X-kX^{k-1}+a^2(m+q)+ca+1 \geq R'(1)=-(a^2(m+q-2)-ca+k-1)=0
\]
if $a=\frac{c}{2(m+q-2)}$ and $c=2\sqrt{p-q}$.

\mbox{}

\noindent {\em II. Case $\: (m+q) \leq 2$:} In this case we set
\[
X=f^k \: , \:\:\:  Y=\sqrt{\frac{m+q}{2}}f^{\frac{m-q-2}{2}}f' \: , \:\:\:
\sqrt{\frac{2}{m+q}}X^{\frac{m-q}{2k}} d\tau=d\xi
\]
where $k= \mbox{\rm min} \{ \frac{2-m-q}{2},p-q \}$ if $m+q \neq 2$, or $k=p-q$ if $m+q=2$, and we arrive at:
\[
\left\{ \begin{array}{ll}
X'=\gamma XY\\ [4pt] Y'=-Y(Y+c_1X^{k_1})+1-X^{k_2}
\end{array}
\right.
\]
with $k_1=1$ and $k_2=\frac{2(p-q)}{2-m-q}$ if $k=\frac{2-m-q}{2}$, or $k_2=1$ and
$k_1=\frac{2-m-q}{2(p-q)}$ if $k=p-q$, where $c_1=c\sqrt{\frac2{m+q}}$ and
$\gamma =\frac{2k}{m+q}$. If $m+q=2$, then $k_1=0$ and $k_2=1$.

\noindent After this point, the proof is entirely analogous to the one of the case $m+q>2$. The demonstration of Theorem 1 is thus complete. \hfill {\sc Q.E.D.}

\mbox{}

\noindent{\sc Remark 2.} {\sl It is worth emphasizing the oscillatory property present in some of the travelling wave solutions identified in Theorem 1, which is to our knowledge a new feature in the context of the solutions of equation (\ref{ec1}). Note also that it is not possible to determine such kind of solutions with the method used by Gilding and Kersner
[\ref{GK}].}

\mbox{}

\noindent{\sc Remark 3.} {\sl Theorem 1 provides also the exact value of $c^*$. Recall that its determination has remained as an open issue since the lower and upper estimations provided by Biro [\ref{B}].}

\mbox{}

\noindent{\sc Remark 4.} {\sl A similar analysis to the one presented in [\ref{PS}] shows that under the same hypotheses of Theorem 1, if $q<1$ and $m>q$ then the travelling waves have the property of finite propagation, namely they vanish for $\xi \geq \xi_0$ for some $\xi _0 < \infty$.}

\pagebreak
\begin{center}
{\bf REFERENCES}
\end{center}
\begin{enumerate}
\item R.A. Fisher, The wave of advance of advantageous genes, {\em Ann. Eugenics\/} {\bf 7}
    355-369 (1937).\label{Fi}
\item A. Kolmogoroff, I. Petrovsky and N. Piscounoff, \'Etude de l'\'equation de la diffusion   avec croissance de la quantit\'e de mati\`ere et son application \`a un probleme
    biologique, {\em Bull. Univ. Moscou, Ser. Internat. Sec. Math.\/} {\bf 1} 1-25
    (1937).\label{KPP}
\item D.G. Aronson and H.F. Weinberger, Nonlinear diffusion in population genetics, combustion,
    and nerve pulse propagation, In {\em Partial Differential Equations and Related
    Topics,\/} (Lecture Notes in Mathematics Vol. 446), pp. 5-49, Springer-Verlag, Berlin,
    (1975).\label{AW}
\item B.H. Gilding and R. Kersner, {\em Travelling waves in nonlinear
    diffusion-convection-reaction\/} (Memorandum No. 1585), Faculty of Mathematical Sciences,
    University of Twente, (2001).\label{GK}
\addtolength{\itemsep}{-2mm}
\item Z. Biro, Attractors in a density-dependent diffusion-reaction model, {\em Nonlinear
    Anal.\/} {\bf 29} 485-499 (1997).\label{B}
\item A. de Pablo and J.L. V\'{a}zquez, Travelling waves and finite propagation in a
    reaction-diffusion equation {\em J. Differential Equations\/} {\bf 93} 19-61
    (1991).\label{PV}
\item L. Malaguti and C. Marcelli, Sharp profiles in degenerate and doubly degenerate
    Fisher-KPP equations, {\em J. Differential Equations\/} {\bf 195} 471-496
    (2003).\label{MM1}
\item F. S\'{a}nchez-Gardu\~{n}o and P.K. Maini, Travelling wave phenomena in some degenerate
    reaction-diffusion equations, {\em J. Differential Equations\/} {\bf 117} 281-319
    (1995).\label{SM}
\item E.O. Voit, {\em Computational Analysis of Biochemical Systems: A Practical Guide for  Biochemists and Molecular Biologists,\/} Cambridge University Press, Cambridge UK
    (2000).\label{VO}
\item N.V. Torres and E.O. Voit, {\em Pathway Analysis and Optimization in Metabolic
    Engineering,\/} Cambridge University Press, Cambridge UK (2002).\label{NV}
\item J.D. Murray, {\em Mathematical Biology,\/} Springer-Verlag, New York (1993).\label{M}
\item A. de Pablo and A. S\'{a}nchez, Global travelling waves in reaction-convection-diffusion
    equations, {\em J. Differential Equations\/} {\bf 165} 377-413 (2000).\label{PS}
\end{enumerate}
\end{document}